\documentclass[12pt]{amsart}

\usepackage[utf8]{inputenc}
\usepackage{amssymb}
\usepackage{amsfonts}
\usepackage{amsmath}
\usepackage{amsthm}
\usepackage{tikz}
\usepackage{tikz-cd}
\usepackage{xcolor}
\usepackage{graphicx}
\usepackage{array}
\usepackage{mathrsfs}
\usepackage{hyperref}

\theoremstyle{definition}

\begin{document}

\title{Numerical Integrators for Mechanical Systems on Lie Groups}

\author{Viyom Vivek}
\address{Centre for Systems and Control, IIT Bombay, India}
\email{viyomvivek@iitb.ac.in}
\author{David Martin de Diego}
\address{ICMAT Madrid, Spain}
\email{david.martin@icmat.es}

\author{Ravi Banavar}
\address{Centre for Systems and Control, IIT Bombay, India}
\email{banavar@iitb.ac.in}

\begin{abstract}
Retraction maps are known to be the seed for all numerical integrators. These retraction maps-based integrators can be further lifted to tangent and cotangent bundles, giving rise to structure-preserving integrators for mechanical systems. We explore the particular case where the configuration space of our mechanical system is a Lie group with certain symmetries. Here, the integrator simplifies based on the property that the tangent and cotangent bundles of Lie groups are trivializable. Finally, we present a framework for designing numerical integrators for Euler-Poincar{\'e} and Lie-Poisson type equations.
\end{abstract}

\maketitle

\tableofcontents

\section{Trivializations on Lie Groups}

\subsection{Trivialization of the tangent and cotangent bundle of a Lie group}

Let $G$ be a Lie group with the Lie algebra denoted by $Lie(G)$ or $\mathfrak{g}$. The group structure on $G$ induces a canonical group structure on the tangent bundle of $G$ denoted by $TG$, see \cite{c3} or \cite{c4}. This induced group multiplication and inversion is given by
\begin{align*}
    u_g \, v_h &:= (T_gR_h(u)+T_hL_g(v))_{gh} = (u \cdot h+g \cdot v)_{gh} \quad \text{and} \\
    (u_g)^{-1} &:= (-T_eR_{g^{-1}} \circ T_gL_{g^{-1}}(u))_{g^{-1}} = (-g^{-1} \cdot u \cdot g^{-1})_{g^{-1}}
\end{align*}
respectively. Here $u_g, v_h \in TG$ and we will use the shorthand notation interchangeably. It is known that the tangent bundle of any Lie group can be trivialized using left (or right) translation as
\begin{align*}
    tr_{TG} : TG &\to G \times \mathfrak{g} \\
    u_g &\mapsto (g, \xi)
\end{align*}
where $\xi:=T_gL_{g^{-1}}(u)=g^{-1} \cdot u$. We would like to equip $G \times \mathfrak{g}$ with a group structure such that $tr_{TG}$ becomes a group isomorphism on top of being a diffeomorphism. So,
\begin{align*}
    &tr_{TG}(u_g \, v_h) = tr_{TG}(u_g)\, tr_{TG}(v_h) \\
    &\implies tr_{TG}((u \cdot h+g \cdot v)_{gh}) = (g,g^{-1}\cdot u) (h, h^{-1}\cdot v) \\
    &\implies (gh,(gh)^{-1} \cdot (u \cdot h + g \cdot v)) = (g,g^{-1}\cdot u) (h, h^{-1}\cdot v) \\
    &\implies (gh,h^{-1}g^{-1}\cdot u \cdot h+ h^{-1}\cdot v)=(g,g^{-1}\cdot u) (h, h^{-1}\cdot v).
\end{align*}
Hence, the group multiplication on $G \times \mathfrak{g}$ is given by 
\begin{align*}
    (g,\xi)(h,\eta):=(gh,Ad_{h^{-1}}\xi+\eta)
\end{align*}
where $(g,\xi),(h,\eta) \in G \times \mathfrak{g}$. The identity element is the zero vector of the Lie algebra i.e. $(e,0)$ while the inverse is given by
\begin{align*}
    (g,\xi)^{-1}:=(g^{-1},-Ad_g\xi).
\end{align*}
Notice that this defines a semidirect product structure on $G\times \mathfrak{g}$ where $G$ is treated as a group while $\mathfrak{g}$ is treated as a vector space with $G$ acting through the adjoint action. Henceforth we will denote this group as $G \ltimes_L \mathfrak{g}$ or simply $G \ltimes \mathfrak{g}$. 

Analogously, the cotangent bundle of $G$ denoted by $T^*G$ can be first equipped with the induced canonical group structure and then trivialized using left (or right) translation. We end up with $G \ltimes_L \mathfrak{g}^*$ or $G \ltimes \mathfrak{g}^*$ where the group multiplication is given by
\begin{align*}
    (g,\mu)(h,\nu):=(gh,Ad^*_{h^{-1}}\mu+\nu)
\end{align*}
where $(g,\mu),(h,\nu) \in G \times \mathfrak{g}^*$. Here the identity element is the zero vector of the dual of the Lie algebra i.e. $(e,0)$ while the inverse is given by
\begin{align*}
    (g,\mu)^{-1}:=(g^{-1},-Ad^*_g\mu).
\end{align*}

\subsection{Trivialization of higher order bundles of a Lie group}
The higher order bundles of a Lie group $G$ namely $TTG,TT^*G,T^*TG$ and $T^*T^*G$ can be trivialized in multiple ways, see \cite{c1}, \cite{c5} or \cite{c6}. Let us begin with trivializing $TTG$. We first trivialize $TG$ to get the following identification $TTG \cong T(G \ltimes \mathfrak{g})$ and then trivialize $T(G \ltimes \mathfrak{g}) \cong (G \ltimes \mathfrak{g})\ltimes Lie(G \ltimes \mathfrak{g})$ or $(G \ltimes \mathfrak{g})\ltimes (\mathfrak{g} \ltimes \mathfrak{g})$ to finally obtain $TTG \cong (G \ltimes \mathfrak{g})\ltimes (\mathfrak{g} \ltimes \mathfrak{g})$. The group multiplication is given by
\begin{align*}
    (g,\xi,\Bar{\xi},\Tilde{\xi})(h,\eta,\Bar{\eta},\Tilde{\eta}):=&(gh,Ad_{h^{-1}}\xi+\eta, \\
    &Ad_{h^{-1}}\Bar{\xi}+\Bar{\eta},Ad_{h^{-1}}(\Tilde{\xi}+[\Bar{\xi},Ad_h\eta])+\Tilde{\eta})
\end{align*}
where $(g,\xi,\Bar{\xi},\Tilde{\xi}),(h,\eta,\Bar{\eta},\Tilde{\eta}) \in (G \ltimes \mathfrak{g})\ltimes (\mathfrak{g} \ltimes \mathfrak{g})$. The identity is $(e,0,0,0)$ and the inverse is given by
\begin{align*}
    (g,\xi,\Bar{\xi},\Tilde{\xi})^{-1}:=(g^{-1},-Ad_g\xi,-Ad_g\Bar{\xi},-Ad_g(\Tilde{\xi}+[\xi,\Bar{\xi}])).
\end{align*}
Note that we will be using left trivialization throughout without stating it explicitly. In a similar fashion, we get $TT^*G \cong (G \ltimes \mathfrak{g}^*) \ltimes (\mathfrak{g} \ltimes \mathfrak{g}^*)$, $T^*TG \cong (G \ltimes \mathfrak{g}) \ltimes (\mathfrak{g}^* \ltimes \mathfrak{g}^*)$ and $T^*T^*G \cong (G \ltimes \mathfrak{g}^*) \ltimes (\mathfrak{g}^* \ltimes \mathfrak{g})$.

\subsection{Trivialized Tulczyjew's triple}
Recall that on any manifold $M$ we can define three bundle morphisms namely $\alpha_M: TT^*M \to T^*TM$, $\beta_M: TT^*M \to T^*T^*M$ and $\kappa_M: TTM \to TTM$ where the first two also happen to be symplectomorhisms when equipped with canonical symplectic structures, see \cite{c9}, \cite{c10}, \cite{c11} or \cite{c8}. These can be summarized in a couple of commutative diagrams as shown below.
\begin{center}
\begin{tikzcd}
T^*T^*M \arrow[dr,"\pi_{T^*M}",swap] & & TT^*M \arrow[ll,"\beta_M",swap] \arrow[rr,"\alpha_M"] \arrow[dl,"\tau_{T^*M}"] \arrow[dr,"T\pi_M",swap] & & T^*TM \arrow[dl,"\pi_{TM}"] \\
& T^*M \arrow[dr,"\pi_M",swap] & & TM \arrow[dl,"\tau_M"] & \\
& & M & &
\end{tikzcd}
\begin{tikzcd}
TTM \arrow[rr,"\kappa_M",shift left] \arrow[dr,"\tau_{TM}",swap] & & TTM \arrow[ll,"\kappa_M",shift left] \arrow[dl,"T\tau_M"] \\
& TM \arrow[d,"\tau_M"] & \\
& M &
\end{tikzcd}
\end{center}
In the particular case of Lie groups i.e. $M=G$, we can obtain a trivialized version of the Tulczyjew's triple given by
\begin{align*}
    \Tilde{\alpha}_G(g,\mu,\xi,\nu)&:=(g,\xi,\nu+ad_\xi^*\mu,\mu),\\
    \Tilde{\beta}_G(g,\mu,\eta,\nu)&:=(g,\mu,\nu+ad_\eta^*\mu,-\eta) \quad \text{and}\\
    \Tilde{\kappa}_G(g,\xi,\eta,\zeta)&:=(g,\eta,\xi,\zeta+[\xi,\eta])
\end{align*}
where $g \in G$, $\xi,\eta,\zeta \in \mathfrak{g}$ and $\mu,\nu \in \mathfrak{g}^*$. These too can be summarized in the commutative diagrams shown below.
\begin{center}
    \begin{tikzpicture}[scale=0.9]
% grid
%\draw[dotted] (-7,-6) grid (7,4);
\node at (0,0) {$(G \ltimes \mathfrak{g}^*) \ltimes(\mathfrak{g} \ltimes \mathfrak{g}^*)$};
\draw[->] (2,0)--(3,0);
\node at (2.5,0.5) {$\Tilde{\alpha}_G$};
\draw[->] (-2,0)--(-3,0);
\node at (-2.5,0.5) {$\Tilde{\beta}_G$};
\draw[->] (-0.5,-0.5)--(-2,-1.5);
\node at (-1.8,-0.8) {$pr_{12}$};
\draw[->] (0.5,-0.5)--(2,-1.5);
\node at (1.8,-0.8) {$pr_{13}$};
\node at (-5,0) {$(G \ltimes \mathfrak{g}^*) \ltimes(\mathfrak{g}^* \ltimes \mathfrak{g})$};
\draw[->] (-4.5,-0.5)--(-3,-1.5);
\node at (-4.2,-1.2) {$pr_{12}$};
\draw[->] (-2,-2.5)--(-0.5,-3.5);
\node at (-1.5,-3.2) {$pr_1$};
\draw[->] (4.5,-0.5)--(3,-1.5);
\node at (4.2,-1.2) {$pr_{12}$};
\draw[->] (2,-2.5)--(0.5,-3.5);
\node at (1.5,-3.2) {$pr_1$};
\node at (5,0) {$(G \ltimes \mathfrak{g}) \ltimes(\mathfrak{g}^* \ltimes \mathfrak{g}^*)$};
\node at (-2.5,-2) {$G \ltimes \mathfrak{g}^*$};
%\draw[blue,->] (-3,-2.5)--(-4.5,-3.5);
%\node at (-3.5,-3.2) {\textcolor{blue}{$h$}};
%\node at (-5,-4) {$\mathbb{R}$};
\node at (2.5,-2) {$G \ltimes \mathfrak{g}$};
\node at (0,-4) {$G$};
%\draw[red,->] (3,-2.5)--(4.5,-3.5);
%\node at (3.5,-3.2) {\textcolor{red}{$l$}};
%\node at (5,-4) {$\mathbb{R}$};

%\node at (5,1) {\textcolor{red}{$(g,\xi,ad_\xi^*(\frac{\delta l}{\delta \xi}),\frac{\delta l}{\delta \xi})$}};
%\node at (-5,1) {\textcolor{blue}{$(g,\mu,ad_{\frac{\delta h}{\delta \mu}}^*(\mu),\frac{\delta h}{\delta \mu})$}};
%\node at (5,1) {\textcolor{red}{$(g,\xi,ad_\xi^*(\frac{\delta l}{\delta \xi}),\frac{\delta l}{\delta \xi})$}};
%\node at (0,1) {\textcolor{red}{$(g,\frac{\delta l}{\delta \xi},\xi,0)$}};
%\node at (0,2) {\textcolor{blue}{$(g,\mu,\frac{\delta h}{\delta \mu},0)$}};

%\node at (4,-2) {\textcolor{red}{$(g,\xi)$}};
%\node at (-4,-2) {\textcolor{blue}{$(g,\mu)$}};

%\node at (5,-5) {\textcolor{red}{$l(\xi)$}};
%\node at (-5,-5) {\textcolor{blue}{$h(\mu)$}};
%\node at (0,-5) {$g$};

\end{tikzpicture}
\end{center}
\begin{center}
\begin{tikzcd}[column sep=tiny]
(G \ltimes \mathfrak{g})\ltimes(\mathfrak{g} \ltimes \mathfrak{g}) \arrow[rr,"\Tilde{\kappa}_G",shift left] \arrow[dr,"pr_{12}",swap] & & (G \ltimes \mathfrak{g})\ltimes(\mathfrak{g} \ltimes \mathfrak{g}) \arrow[ll,"\Tilde{\kappa}_G",shift left] \arrow[dl,"pr_{13}"] \\
& G \ltimes \mathfrak{g} \arrow[d,"pr_1"] & \\
& G &
\end{tikzcd}
\end{center}
Here $pr_i$ and $pr_{ij}$ are projections onto $i^{th}$ and $ij^{th}$ factors, respectively.

\section{Mechanical Systems on Lie Groups}

\subsection{Euler-Poincar{\'e} equations}
Let the Lie group $G$ be the configuration space of a mechanical system with the Lagrangian $\mathscr{L}: TG \to \mathbb{R}$. Then the Euler-Lagrange equations are given by 
\begin{align}
    \frac{d}{dt}\left(\frac{\delta \mathscr{L}}{\delta \Dot{g}}\right)=\frac{\delta \mathscr{L}}{\delta g},
\end{align}
see \cite{c14} or \cite{c13}. The trivialized Euler-Lagrange equations are obtained from trivializing $TG$ as
\begin{align}
    \frac{d}{dt}\left(\frac{\delta \Tilde{\mathscr{L}}}{\delta \xi} \right)&= ad^*_\xi\left(\frac{\delta \Tilde{\mathscr{L}}}{\delta \xi}\right)+T^*_eL_g\left(\frac{\delta \Tilde{\mathscr{L}}}{\delta g}\right) \quad \text{and} \\
    \Dot{g}&=T_eL_g(\xi)
\end{align}
where $\Tilde{\mathscr{L}}: G \ltimes \mathfrak{g} \to \mathbb{R}$ is the trivialized Lagrangian and equation (3) is used for reconstruction. In the special case where the trivialized Lagrangian just depends on the Lie algebra i.e. the Lagrangian is left (or right) invariant, the trivialized Euler-Lagrange equations reduce to the Euler-Poincar{\'e} equations given by
\begin{align}
    \frac{d}{dt}\left(\frac{\delta l}{\delta \xi} \right)&= ad^*_\xi\left(\frac{\delta l}{\delta \xi}\right) \quad \text{and}\\
    \Dot{g}&=T_eL_g(\xi)
\end{align}
where $l: \mathfrak{g} \to \mathbb{R}$ is the reduced Lagrangian.

\subsection{Lie-Poisson equations}
Let $\mathfrak{g}^*$ be the dual of the Lie algebra of $G$ equipped with the Lie-Poisson bracket. Then the Lie-Poisson equations corresponding to a reduced Hamiltonian $h: \mathfrak{g}^* \to \mathbb{R}$ are given by
\begin{align}
    \Dot{\mu}&=ad_{\frac{\delta h}{\delta \mu}}^*(\mu) \quad \text{and}\\
    \Dot{g}&=T_eL_g\left(\frac{\delta h}{\delta \mu}\right)
\end{align}
see \cite{c13}, \cite{c14} or \cite{c17}. Here again equation (7) is the reconstruction equation. In the presence of a Legendre transform $\mathbb{F}l:\mathfrak{g} \to \mathfrak{g}^*$ that takes $\xi \mapsto \frac{\delta l}{\delta \xi}=:\mu$ and $h(\mu)=\left<\mu,\xi\right>-l(\xi)$, the set of equations (4) and (5) can be shown to be equivalent to (6) and (7).

\subsection{Euler-Arnold equations}
Let the Lie algebra $\mathfrak{g}$ be equipped with an inner product $\langle\langle \cdot,\cdot \rangle\rangle$ with the corresponding inertia operator $\mathbb{I}: \mathfrak{g} \to \mathfrak{g}^*$ defined by
\begin{align}
    \langle \mathbb{I}(\xi),\eta \rangle:=\langle\langle \xi,\eta \rangle\rangle, 
\end{align}
see \cite{c16}, \cite{c15} or \cite{c17}. Then the set of equations (4) and (5) for the reduced lagrangian $l(\xi):=\frac{1}{2}\langle \mathbb{I}(\xi),\xi \rangle$ become
\begin{align}
    \Dot{\xi}&=\mathbb{I}^{-1}(ad_\xi^*(\mathbb{I}(\xi))) \quad \text{and}\\
    \Dot{g}&=T_eL_g(\xi)
\end{align}
which are called the Euler-Arnold equations. These describe the geodesics on $G$ with respect to the left-invariant metric formed by $\langle\langle \cdot,\cdot \rangle\rangle$.

The Euler-Poincar{\'e} and the Lie-Poisson equations can be interpreted as Lagrangian submanifolds of $TT^*G$ after trivializaion to $(G \ltimes \mathfrak{g}^*)\ltimes(\mathfrak{g} \ltimes \mathfrak{g}^*)$ as shown in the diagram below.
\begin{center}
    \begin{tikzpicture}[scale=0.9]
% grid
%\draw[dotted] (-7,-6) grid (7,4);
%
\node at (0,0) {$(G \ltimes \mathfrak{g}^*) \ltimes(\mathfrak{g} \ltimes \mathfrak{g}^*)$};
\draw[->] (2,0)--(3,0);
\node at (2.5,0.5) {$\Tilde{\alpha}_G$};
\draw[->] (-2,0)--(-3,0);
\node at (-2.5,0.5) {$\Tilde{\beta}_G$};
\draw[->] (-0.5,-0.5)--(-2,-1.5);
\node at (-1.8,-0.8) {$pr_{12}$};
\draw[->] (0.5,-0.5)--(2,-1.5);
\node at (1.8,-0.8) {$pr_{13}$};
\node at (-5,0) {$(G \ltimes \mathfrak{g}^*) \ltimes(\mathfrak{g}^* \ltimes \mathfrak{g})$};
\draw[->] (-4.5,-0.5)--(-3,-1.5);
\node at (-3.2,-0.8) {$pr_{12}$};
\draw[->] (-2,-2.5)--(-0.5,-3.5);
\node at (-1.5,-3.2) {$pr_1$};
\draw[->] (4.5,-0.5)--(3,-1.5);
\node at (3.2,-0.8) {$pr_{12}$};
\draw[->] (2,-2.5)--(0.5,-3.5);
\node at (1.5,-3.2) {$pr_1$};
\node at (5,0) {$(G \ltimes \mathfrak{g}) \ltimes(\mathfrak{g}^* \ltimes \mathfrak{g}^*)$};
\node at (-2.5,-2) {$G \ltimes \mathfrak{g}^*$};
\draw[blue,->] (-3,-2.5)--(-4.5,-3.5);
\node at (-3.5,-3.2) {\textcolor{blue}{$h$}};
\node at (-5,-4) {$\mathbb{R}$};
\node at (2.5,-2) {$G \ltimes \mathfrak{g}$};
\node at (0,-4) {$G$};
\draw[red,->] (3,-2.5)--(4.5,-3.5);
\node at (3.5,-3.2) {\textcolor{red}{$l$}};
\node at (5,-4) {$\mathbb{R}$};

\node at (5,1) {\textcolor{red}{$(g,\xi,ad_\xi^*(\frac{\delta l}{\delta \xi}),\frac{\delta l}{\delta \xi})$}};
\node at (-5,1) {\textcolor{blue}{$(g,\mu,ad_{\frac{\delta h}{\delta \mu}}^*(\mu),\frac{\delta h}{\delta \mu})$}};
\node at (5,1) {\textcolor{red}{$(g,\xi,ad_\xi^*(\frac{\delta l}{\delta \xi}),\frac{\delta l}{\delta \xi})$}};
\node at (0,1) {\textcolor{red}{$(g,\frac{\delta l}{\delta \xi},\xi,0)$}};
\node at (0,2) {\textcolor{blue}{$(g,\mu,\frac{\delta h}{\delta \mu},0)$}};

\node at (4,-2) {\textcolor{red}{$(g,\xi)$}};
\node at (-4,-2) {\textcolor{blue}{$(g,\mu)$}};

\node at (5,-5) {\textcolor{red}{$l(\xi)$}};
\node at (-5,-5) {\textcolor{blue}{$h(\mu)$}};
\node at (0,-5) {$g$};

\end{tikzpicture}
\end{center}

\section{Discretization on Lie Groups}
\subsection{Retraction Maps}
Let $M$ be a smooth manifold with its tangent bundle denoted by $TM$. Then, a smooth map $R:TM \to M$ satisfying (i) $R(0_x)=x$ and (ii) $\left.\frac{d}{dt}\right|_{t=0}R(tv_x)=v_x$ for every $x \in M$ and $v_x \in TM$ is called a retraction map on $M$, see \cite{c2} or \cite{c18}. 

Some common examples include (i) $R(v_x):=\exp_\mathcal{G}{(v_x)}$ where $\exp_\mathcal{G}$ denotes the Riemannian exponential on a Riemannian manifold $(M,\mathcal{G})$ and (ii) $R(v_g):=\exp{(v_g)}$ where $\exp$ denotes the group exponential on a Lie group $G$. An extended list of examples of retraction maps can be found in \cite{c18}.
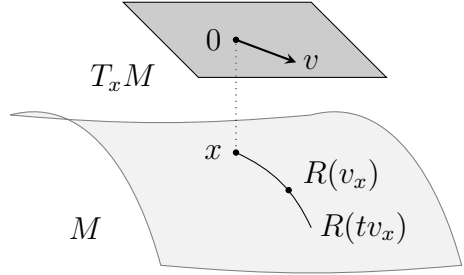
\begin{center}
    \begin{tikzpicture}[scale=1]
        % grid
        %\draw[dotted,white] (0,0) grid (8,5);
        % manifold
        \draw[gray,fill=gray,fill opacity=0.1] (1,3) to[out=15,in=105] (3,1) to[out=-5,in=-175] (7,1) to[out=105,in=15] (5,3) to[out=-175,in=-5] (1,3);
        \node at (2,1.5) {$M$};
        % point and curve
        \draw[fill] (4,2.5) circle (1pt);
        \draw (4,2.5) to[out=-25,in=115] (5,1.5);
        \draw[fill] (4.7,2) circle (1pt);
        \node at (3.7,2.5) {$x$};
        \node at (5.7,1.5) {$R(tv_x)$};
        \node at (5.4,2.2) {$R(v_x)$};
        % tangent space
        \draw[dotted] (4,2.5)--(4,4);
        \draw[fill=gray,fill opacity=0.4] (3.5,3.5)--(6,3.5)--(5,4.5)--(2.5,4.5)--cycle;
        \draw[fill] (4,4) circle (1pt);
        \draw[thick,-stealth] (4,4)--(4.8,3.7);
        \node at (3.7,4) {$0$};
        \node at (5,3.7) {$v$};
        \node at (2.5,3.5) {$T_xM$};
    \end{tikzpicture}
\end{center}

We will use (left) trivialized retraction maps on a Lie group $G$ defined as $\Tilde{R}: G \ltimes \mathfrak{g} \to G$ which maps $(g,\xi)\mapsto g\tau(\xi)$. Here $\tau: \mathfrak{g} \to G$ is a local diffeomorphism between $0 \in \mathfrak{g}$ and $e \in G$ satisfying $\tau(\xi)\tau(-\xi)=e$ for every $\xi \in \mathfrak{g}$, see \cite{c19} for more details. Right trivialized retraction maps can be defined analogously.
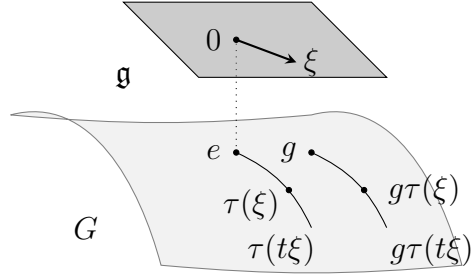
\begin{center}
    \begin{tikzpicture}[scale=1]
        % grid
        %\draw[dotted,white] (0,0) grid (8,5);
        % manifold
        \draw[gray,fill=gray,fill opacity=0.1] (1,3) to[out=15,in=105] (3,1) to[out=-5,in=-175] (7,1) to[out=105,in=15] (5,3) to[out=-175,in=-5] (1,3);
        \node at (2,1.5) {$G$};
        % point and curve
        \draw[fill] (4,2.5) circle (1pt);
        \draw (4,2.5) to[out=-25,in=115] (5,1.5);
        \draw[fill] (4.7,2) circle (1pt);
        \node at (3.7,2.5) {$e$};
        \node at (4.6,1.2) {$\tau(t\xi)$};
        \node at (4.2,1.8) {$\tau(\xi)$};
        % g and curve
        \draw[fill] (5,2.5) circle (1pt);
        \draw (5,2.5) to[out=-25,in=115] (6,1.5);
        \draw[fill] (5.7,2) circle (1pt);
        \node at (4.7,2.5) {$g$};
        \node at (6.6,1.2) {$g\tau(t\xi)$};
        \node at (6.5,2) {$g\tau(\xi)$};
        % tangent space
        \draw[dotted] (4,2.5)--(4,4);
        \draw[fill=gray,fill opacity=0.4] (3.5,3.5)--(6,3.5)--(5,4.5)--(2.5,4.5)--cycle;
        \draw[fill] (4,4) circle (1pt);
        \draw[thick,-stealth] (4,4)--(4.8,3.7);
        \node at (3.7,4) {$0$};
        \node at (5,3.7) {$\xi$};
        \node at (2.5,3.5) {$\mathfrak{g}$};
    \end{tikzpicture}
\end{center}

The most common choice for the $\tau$ map is group exponential but other interesting choices also exist for specific groups which we shall see later. 

\subsection{Discretization Maps}
Let $M$ be a smooth manifold with its tangent bundle denoted by $TM$. Then, a smooth map $R_d:TM \to M \times M$ satisfying (i) $R_d(0_x)=(R^1(0_x),R^2(0_x))=(x,x)$ and (ii) $\left.\frac{d}{dt}\right|_{t=0}R^2(tv_x)-\left.\frac{d}{dt}\right|_{t=0}R^1(tv_x)=v_x$ for every $x \in M$ and $v_x \in TM$ is called a discretization map on $M$, see \cite{c2}.
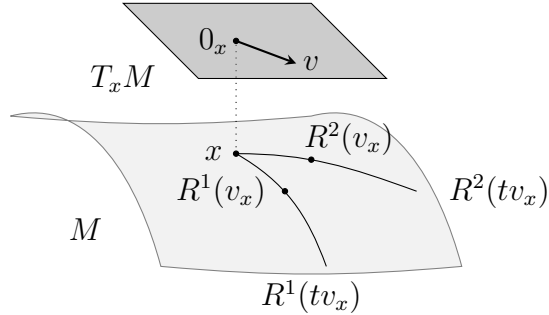
\begin{center}
    \begin{tikzpicture}[scale=1]
        % grid
        %\draw[dotted,white] (0,0) grid (9,5);
        % manifold
        \draw[gray,fill=gray,fill opacity=0.1] (1,3) to[out=15,in=105] (3,1) to[out=-5,in=-175] (7,1) to[out=105,in=15] (5,3) to[out=-175,in=-5] (1,3);
        \node at (2,1.5) {$M$};
        % point and curve
        \draw[fill] (4,2.5) circle (1pt);
        \draw[] (4,2.5) to[out=-30,in=110] (5.2,1);
        \draw[] (4,2.5) to[out=0,in=160] (6.4,2);
        %\draw (4,2.5) to[out=-25,in=115] (5,1.5);
        \draw[fill] (4.65,2) circle (1pt);
        \draw[fill] (5,2.42) circle (1pt);
        \node at (3.7,2.5) {$x$};
        \node at (5.5,2.7) {$R^2(v_x)$};
        \node at (3.8,2) {$R^1(v_x)$};
        \node at (5,0.6) {$R^1(tv_x)$};
        \node at (7.5,2) {$R^2(tv_x)$};
        % tangent space
        \draw[dotted] (4,2.5)--(4,4);
        \draw[fill=gray,fill opacity=0.4] (3.5,3.5)--(6,3.5)--(5,4.5)--(2.5,4.5)--cycle;
        \draw[fill] (4,4) circle (1pt);
        \draw[thick,-stealth] (4,4)--(4.8,3.7);
        \node at (3.7,4) {$0_x$};
        \node at (5,3.7) {$v$};
        \node at (2.5,3.5) {$T_xM$};
    \end{tikzpicture}
\end{center}

Note that these two properties ensure that $R_d:TM \to M \times M$ is a local diffeomorphism. Also notice that given a retraction map $R: TM \to M$, one can construct a family of discretization maps $R_d:TM \to M \times M$ as $R_d(v_x):=(R(-\theta v_x),R((1-\theta)v_x))$ where $\theta \in [0,1]$ is some parameter. The proofs can be found in \cite{c2}. 
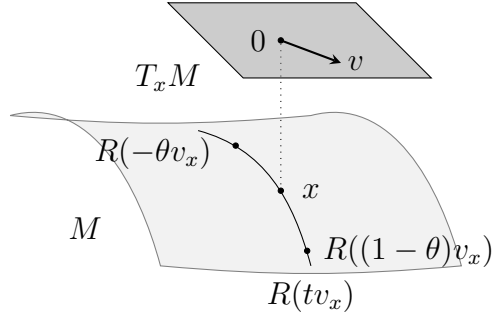
\begin{center}
    \begin{tikzpicture}
        % grid
        %\draw[dotted,white] (0,0) grid (8,5);
        % manifold
        \draw[gray,fill=gray,fill opacity=0.1] (1,3) to[out=15,in=105] (3,1) to[out=-5,in=-175] (7,1) to[out=105,in=15] (5,3) to[out=-175,in=-5] (1,3);
        \node at (2,1.5) {$M$};
        % point and curve
        \draw[fill] (4,2.6) circle (1pt);
        \draw (3.5,2.8) to[out=-15,in=105] (5,1);
        \draw[fill] (4.6,2) circle (1pt);
        \draw[fill] (4.95,1.2) circle (1pt);
        \node at (2.9,2.5) {$R(-\theta v_x)$};
        \node at (6.3,1.2) {$R((1-\theta)v_x)$};
        \node at (5,2) {$x$};
        % tangent space
        \draw[dotted] (4.6,2)--(4.6,4);
        \draw[fill=gray,fill opacity=0.4] (4.1,3.5)--(6.6,3.5)--(5.6,4.5)--(3.1,4.5)--cycle;
        \draw[fill] (4.6,4) circle (1pt);
        \draw[thick,-stealth] (4.6,4)--(5.4,3.7);
        \node at (4.3,4) {$0$};
        \node at (5.6,3.7) {$v$};
        \node at (3.1,3.5) {$T_xM$};
        \node at (5,0.6) {$R(tv_x)$};
    \end{tikzpicture}
\end{center}
Again, we will use (left) trivialized discretization maps on a Lie group $G$ defined as $\Tilde{R}_d:G \ltimes \mathfrak{g} \to G \times G$ which maps $(g,\xi) \mapsto (g\tau(-\theta \xi),g\tau((1-\theta)\xi))$ for some $\theta \in [0,1]$.
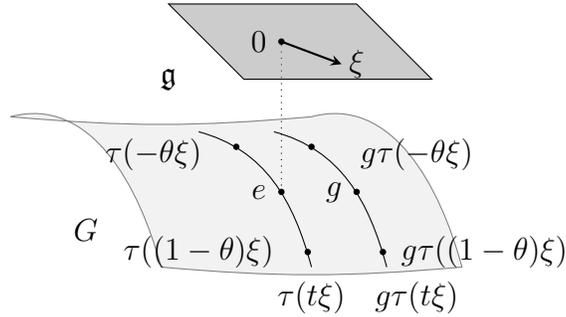
\begin{center}
    \begin{tikzpicture}
        % grid
        %\draw[dotted,white] (0,0) grid (8,5);
        % manifold
        \draw[gray,fill=gray,fill opacity=0.1] (1,3) to[out=15,in=105] (3,1) to[out=-5,in=-175] (7,1) to[out=105,in=15] (5,3) to[out=-175,in=-5] (1,3);
        \node at (2,1.5) {$G$};
        % point and curve
        \draw[fill] (4,2.6) circle (1pt);
        \draw (3.5,2.8) to[out=-15,in=105] (5,1);
        \draw[fill] (4.6,2) circle (1pt);
        \draw[fill] (4.95,1.2) circle (1pt);
        \node at (2.9,2.5) {$\tau(-\theta\xi)$};
        \node at (3.5,1.2) {$\tau((1-\theta)\xi)$};
        \node at (4.3,2) {$e$};
        % g and curve
        \draw[fill] (5,2.6) circle (1pt);
        \draw (4.5,2.8) to[out=-15,in=105] (6,1);
        \draw[fill] (5.6,2) circle (1pt);
        \draw[fill] (5.95,1.2) circle (1pt);
        \node at (6.4,2.5) {$g\tau(-\theta\xi)$};
        \node at (7.3,1.2) {$g\tau((1-\theta)\xi)$};
        \node at (5.3,2) {$g$};
        % tangent space
        \draw[dotted] (4.6,2)--(4.6,4);
        \draw[fill=gray,fill opacity=0.4] (4.1,3.5)--(6.6,3.5)--(5.6,4.5)--(3.1,4.5)--cycle;
        \draw[fill] (4.6,4) circle (1pt);
        \draw[thick,-stealth] (4.6,4)--(5.4,3.7);
        \node at (4.3,4) {$0$};
        \node at (5.6,3.7) {$\xi$};
        \node at (3.1,3.5) {$\mathfrak{g}$};
        \node at (5,0.6) {$\tau(t\xi)$};
        \node at (6.4,0.6) {$g\tau(t\xi)$};
    \end{tikzpicture}
\end{center}

\subsection{Numerical Integrators on a Manifold}
Let $M$ be a smooth manifold and $X:M \to TM$ be a smooth vector field whose dynamics we wish to discretize. We use a discretization map $R_d:TM \to M \times M$ to achieve this as highlighted by the following commutative diagram, see \cite{c2}.
\begin{center}
    \begin{tikzcd}[row sep=large,column sep=large]
        M \times M \arrow[r,"R_d^{-1}"] \arrow[d,swap,"R_d^{-1}"] & TM \arrow[d,"\tau_M"] \\
        TM & M \arrow[l,"X"]
    \end{tikzcd}
\end{center}
\begin{align}
    tX(\tau_M(R_d^{-1}(x_k,x_{k+1})))=R_d^{-1}(x_k,x_{k+1})
\end{align}
Here $\tau_M:TM \to M$ is the canonical tangent bundle projection and $t \in \mathbb{R}$ is the step size. Note that the step size $t$ plays the crucial role of scaling down the vector field $X$ to bring it within the domain of invertibility of $R_d$ and hence ensuring that (11) is well defined. Now, we can iteratively solve the difference equation (11) given the initial conditions $x_0$.

\subsection{Numerical Integrators on Tangent and cotangent Bundles}
As one only requires a discretization map to develop a numerical integrator on any manifold, we tangent lift $R_d: TM \to M \times M$ to $R_d^T:TTM \to TM \times TM$ such that the following diagram commutes.
\begin{center}
\begin{tikzcd}
TTM \arrow[r,"R_d^T"] \arrow[d,"\kappa_M",swap] & TM \times TM \arrow[d,equal] \\
TTM \arrow[d,"\tau_{TM}",swap] \arrow[r,"TR_d"] & T (M \times M) \arrow[d,"\tau_{M\times M}"] \\
TM \arrow[r,"R_d",swap] & M \times M
\end{tikzcd}
\end{center}
\begin{align}
    R_d^T:=TR_d \circ \kappa_M
\end{align}
Here $\kappa_M$ is one of Tulczyjew's bundle morphism which is commonly known as the canonical involution or canonical flip. Now, the flow of a vector field $X:TM \to TTM$ can be discretized as highlighted by the following commutative diagram.
\begin{center}
    \begin{tikzcd}[row sep=large,column sep=large]
        TM \times TM \arrow[r,"(R_d^{T})^{-1}"] \arrow[d,swap,"(R_d^{T})^{-1}"] & TTM \arrow[d,"\tau_{TM}"] \\
        TTM & TM \arrow[l,"X"]
    \end{tikzcd}
\end{center}
\begin{align}
    tX(\tau_{TM}((R_d^{T})^{-1}(x_k,v_k;x_{k+1},v_{k+1})))=(R_d^{T})^{-1}(x_k,v_k;x_{k+1},v_{k+1})
\end{align}

Similarly, we cotangent lift $R_d:TM \to M \times M$ to $R_d^{T^*}:TT^*M \to T^*M \times T^*M$ such that the following diagram commutes.
\begin{center}
\begin{tikzcd}
TT^*M \arrow[r,"R_d^{T^*}"] \arrow[d,"\alpha_M",swap] & T^*M \times T^*M \\
T^*TM \arrow[d,"\pi_{TM}",swap] \arrow[r,"T^*R_d^{-1}"] & T^* (M \times M) \arrow[u,"\Phi",swap] \arrow[d,"\pi_{M\times M}"] \\
TM \arrow[r,"R_d",swap] & M \times M
\end{tikzcd}
\end{center}
\begin{align}
    R_d^{T^*}:=\Phi \circ T^*R_d^{-1} \circ \alpha_M
\end{align}
Here $\alpha_M$ is another one of Tulczyjew's bundle morphism defined as the ``dual'' to $\kappa_M$ while $\Phi$ is a symplectomorphism, see \cite{c2} or \cite{c20} for more details on twisted symplectic structure. Again, the flow of a vector field $X:T^*M \to TT^*M$ can be discretized as highlighted by the following commutative diagram.
\begin{center}
    \begin{tikzcd}[row sep=large,column sep=large]
        T^*M \times T^*M \arrow[r,"(R_d^{T^*})^{-1}"] \arrow[d,swap,"(R_d^{T^*})^{-1}"] & TT^*M \arrow[d,"\tau_{T^*M}"] \\
        TT^*M & T^*M \arrow[l,"X"]
    \end{tikzcd}
\end{center}
\begin{align}
    tX(\tau_{T^*M}((R_d^{T^*})^{-1}(x_k,p_k;x_{k+1},p_{k+1})))=(R_d^{T^*})^{-1}(x_k,p_k;x_{k+1},p_{k+1})
\end{align}

\subsection{Numerical Integrators on Tangent and Cotangent Bundles of Lie Groups}
In the case that our manifold $M$ is a Lie group $G$, we can use the trivialization methods developed thus far to obtain the trivialized versions of $R_d^T$ and $R_d^{T^*}$ respectively as follows.
\begin{center}
\begin{tikzcd}
(G \ltimes \mathfrak{g})\ltimes(\mathfrak{g} \ltimes \mathfrak{g}) \arrow[r,"\Tilde{R}_d^{T}"] \arrow[d,"\Tilde{\kappa}_G",swap] & 
(G\ltimes \mathfrak{g}) \times (G\ltimes \mathfrak{g})  \\
(G \ltimes \mathfrak{g})\ltimes(\mathfrak{g} \ltimes \mathfrak{g}) \arrow[r,"T\Tilde{R_d}"] \arrow[d,"pr_{12}",swap] & 
(G \times G) \ltimes (\mathfrak{g} \times \mathfrak{g}) \arrow[d,"pr_{12}"] \arrow[u,equal] \\
G \ltimes \mathfrak{g} \arrow[r,"\Tilde{R}_d",swap] & G \times G
\end{tikzcd}
\end{center}
\begin{align}
    \Tilde{R}_d^T:=T\Tilde{R}_d \circ \Tilde{\kappa}_G
\end{align}
\begin{center}
\begin{tikzcd}
(G \ltimes \mathfrak{g}^*)\ltimes(\mathfrak{g} \ltimes \mathfrak{g}^*) \arrow[r,"\Tilde{R}_d^{T^*}"] \arrow[d,"\Tilde{\alpha}_G",swap] & 
(G\ltimes \mathfrak{g}^*) \times (G\ltimes \mathfrak{g}^*)  \\
(G \ltimes \mathfrak{g})\ltimes(\mathfrak{g}^* \ltimes \mathfrak{g}^*) \arrow[r,"T^*\Tilde{R_d}^{-1}"] \arrow[d,"pr_{12}",swap] & 
(G \times G) \ltimes (\mathfrak{g}^* \times \mathfrak{g}^*) \arrow[d,"pr_{12}"] \arrow[u,"\Tilde{\Phi}",swap] \\
G \ltimes \mathfrak{g} \arrow[r,"\Tilde{R}_d",swap] & G \times G
\end{tikzcd}
\end{center}
\begin{align}
    \Tilde{R}_d^{T^*}:=\Tilde{\Phi}\circ T^*\Tilde{R}_d^{-1} \circ \Tilde{\alpha}_G
\end{align}

Finally, these trivialized discretization maps can be used for developing numerical integrators for the Euler-Poincar{\'e}, Lie-Poisson and Euler-Arnold equations as follows. Let $X: G \ltimes \mathfrak{g} \to T(G\ltimes \mathfrak{g})$ denote the Euler-Poincar{\'e} vector field which can be trivialized as shown in the diagram below
\begin{center}
    \begin{tikzcd}
         & (G \ltimes \mathfrak{g}) \ltimes (\mathfrak{g} \ltimes \mathfrak{g}) \\
         G \ltimes \mathfrak{g} \arrow[r,"X",swap] \arrow[ur,"\Tilde{X}"] & T(G \ltimes \mathfrak{g}) \arrow[u]
    \end{tikzcd}
\end{center}
then the following diagram defines an integrator.
\begin{center}
\begin{tikzcd}
(G \ltimes \mathfrak{g})\ltimes(\mathfrak{g} \ltimes \mathfrak{g}) & (G\ltimes \mathfrak{g}) \times (G\ltimes \mathfrak{g}) \arrow[l,"(\Tilde{R}_d^{T})^{-1}",swap] \arrow[d,"(\Tilde{R}_d^{T})^{-1}"] \\
G \ltimes \mathfrak{g} \arrow[u,"\Tilde{X}"] & (G \ltimes \mathfrak{g})\ltimes(\mathfrak{g} \ltimes \mathfrak{g}) \arrow[l,"pr_{12}"]
\end{tikzcd}
\end{center}
\begin{align}
    t\tilde{X}(pr_{12}((\Tilde{R}_d^{T})^{-1}(g_k,\xi_k;g_{k+1},\xi_{k+1})))=(\Tilde{R}_d^{T})^{-1}(g_k,\xi_k;g_{k+1},\xi_{k+1})
\end{align}
Here $t \in \mathbb{R}$ is the step size. Similarly, let $X:G \ltimes \mathfrak{g}^*\to T(G\ltimes \mathfrak{g}^*)$ denote the Lie-Poisson vector field which can also be trivialized as shown in the diagram below 
\begin{center}
    \begin{tikzcd}
        & (G \ltimes \mathfrak{g}^*) \ltimes (\mathfrak{g} \ltimes \mathfrak{g}^*) \\
         G \ltimes \mathfrak{g}^* \arrow[r,"X",swap] \arrow[ur,"\Tilde{X}"] & T(G \ltimes \mathfrak{g}^*) \arrow[u]
    \end{tikzcd}
\end{center}
then the following diagram defines an integrator.
\begin{center}
\begin{tikzcd}
(G \ltimes \mathfrak{g}^*)\ltimes(\mathfrak{g} \ltimes \mathfrak{g}^*) & (G\ltimes \mathfrak{g}^*) \times (G\ltimes \mathfrak{g}^*) \arrow[l,"(\Tilde{R}_d^{T^*})^{-1}",swap] \arrow[d,"(\Tilde{R}_d^{T^*})^{-1}"] \\
G \ltimes \mathfrak{g}^* \arrow[u,"\tilde{X}"] & (G \ltimes \mathfrak{g}^*)\ltimes(\mathfrak{g} \ltimes \mathfrak{g}^*) \arrow[l,"pr_{12}"]
\end{tikzcd}
\end{center}
\begin{align}
    t\tilde{X}(pr_{12}((\Tilde{R}_d^{T^*})^{-1}(g_k,\mu_k;g_{k+1},\mu_{k+1}))=(\Tilde{R}_d^{T^*})^{-1}(g_k,\mu_k;g_{k+1},\mu_{k+1})
\end{align}
Here again $t \in \mathbb{R}$ is the step size.

\section{Illustrative Example}
Let us carefully develop a numerical integrator for the Lie-Poisson equations on a Lie group $G$. We begin by choosing a discretization map $\Tilde{R}_d : G \ltimes \mathfrak{g} \to G \times G$ as $\Tilde{R}_d(g,\xi) := (g,g\tau(\xi))$. Then, we evaluate the tangent of $\Tilde{R}_d$ at $(g,\xi)$ to get
\begin{align*}
    T_{(g,\xi)}\tilde{R}_d(v,\eta)=(v,T_{\tau(\xi)}L_g(T_\xi\tau(\eta))+T_gR_{\tau(\xi)}(v))
\end{align*}
where $(v,\eta) \in T_{(g,\xi)}(G \ltimes \mathfrak{g})$. Using the natural pairing between dual vectors, we can immediately find
\begin{align*}
    T_{(g,\xi)}^*\tilde{R}_d(p,q)=(p+T^*_gR_{\tau(\xi)}(q),T^*_\xi\tau(T^*_{\tau(\xi)}L_g(q)))
\end{align*}
where $(p,q) \in T^*_{(g,g\tau(\xi))}(G \times G)$. Now, we trivialize both $T^*(G \ltimes \mathfrak{g})$ and $T^*(G \times G)$ as
\begin{align*}
    &tr_{T^*(G \ltimes  \mathfrak{g})}(g,\xi;p+T^*_gR_{\tau(\xi)}(q),T^*_\xi\tau(T^*L_g(q))) =
    (g,\xi;T^*_eL_g(p)\\
    &+T^*_eL_g(T^*_gR_{\tau(\xi)}(q))
    +ad^*_\xi(T^*_\xi\tau(T^*_{\tau(\xi)}L_g(q)),
    T^*_\xi\tau(T^*_{\tau(\xi)}L_g(q))) \quad \text{and}\\
    &tr_{T^*(G \times G)}(g,g\tau(\xi);p,q) = (g,g\tau(\xi);T^*_eL_g(p),
    T^*_eL_{g\tau(\xi)}(q))
\end{align*}
respectively. Finally, we use the trivialized symplectomorphisms $\tilde{\alpha}_G$ and $\tilde{\Phi}$ to compute $\tilde{R}_d^{T^*}$. The expression for $\tilde{\Phi}$ is given by
\begin{align*}
    \tilde{\Phi}(g,\mu;h,\nu):=(g,h;\mu,-\nu)
\end{align*}
where $g,h \in G$ and $\mu,\nu \in \mathfrak{g}^*$.
As we only need the inverse of $\tilde{R}_d^{T^*}$ for the construction of our integrator, we compute it directly to find
\begin{align*}
    &(\tilde{R}_d^{T^*})^{-1}(g_k,\mu_k;g_{k+1},\mu_{k+1})
    =(g_k,-T^*_{\tau^{-1}(g_k^{-1}g_{k+1})}\tau(T^*_{g_k^{-1}g_{k+1}}L_{g_{k+1}^{-1}g_k}(\mu_{k+1}));\\
    &\tau^{-1}(g_k^{-1}g_{k+1}),\mu_k-Ad^*_{g_{k+1}^{-1}g_k}(\mu_{k+1}))\\
    &=(g_k,-d_{\tau^{-1}(g_k^{-1}g_{k+1})}^{L*}\tau(\mu_{k+1}),\tau^{-1}(g_k^{-1}g_{k+1}),
    \mu_k-Ad^*_{g_{k+1}^{-1}g_k}(\mu_{k+1})).
\end{align*}
Here, we have simplified the notation by introducing the left trivialized tangent of $\tau$ as
\begin{align*}
    d^L_\xi\tau: \mathfrak{g} &\longrightarrow \mathfrak{g}\\
    \eta &\longmapsto T_{\tau(\xi)}L_{\tau(-\xi)}(T_\xi\tau(\eta))
\end{align*}
with its dual denoted by
\begin{align*}
    d^{L*}_{\xi}\tau : \mathfrak{g}^* &\longrightarrow \mathfrak{g}^*\\
    \mu &\longmapsto T^*_{\xi}\tau(T^*_{\tau(\xi)}L_{\tau(-\xi)}(\mu))
\end{align*}
where $\xi,\eta \in \mathfrak{g}$ and $\mu \in \mathfrak{g}^*$.
The right trivialized tangent of $\tau$ and its dual can be defined analogously.
The flow of all intermediate computations is summarized in the diagram below.
\begin{center}
    \begin{tikzcd}
        (G \ltimes \mathfrak{g}^*) \ltimes (\mathfrak{g} \ltimes \mathfrak{g}^*) & (G \ltimes \mathfrak{g}^*) \times (G \ltimes \mathfrak{g}^*) \arrow[l,"(\tilde{R}_d^{T^*})^{-1}",swap] \arrow[d,"\tilde{\Phi}"]\\
        (G \ltimes \mathfrak{g}) \ltimes (\mathfrak{g}^* \times \mathfrak{g}^*) \arrow[u,"\tilde{\alpha}_G^{-1}"] & (G \times G) \ltimes (\mathfrak{g}^* \times \mathfrak{g}^*) \arrow[l,"\widetilde{T^*\tilde{R}_d}",swap]\\
        T^*(G \ltimes \mathfrak{g}^*) \arrow[u,"tr_{T^*(G\ltimes \mathfrak{g})}"] \arrow[d,"\pi_{G \ltimes \mathfrak{g}}",swap] & T^*(G\times G)\arrow[u,"tr_{T^*(G\times G)}",swap] \arrow[l,"T^*\tilde{R}_d",swap] \arrow[d,"\pi_{G\times G}"]\\
        G \ltimes \mathfrak{g} \arrow[r,"\tilde{R}_d",swap] & G \times G
    \end{tikzcd}
\end{center}
The vector field corresponding to the Lie-Poisson equations for the Hamiltonian $h(\mu)=\frac{1}{2}\langle \mu,\mathbb{I}^{-1}(\mu) \rangle$ is given by 
\begin{align*}
    X(g,\mu)=(g,\mu;g \cdot \mathbb{I}^{-1}(\mu),ad^*_{\mathbb{I}^{-1}(\mu)}\mu)
\end{align*}
which can be trivialized as shown below
\begin{center}
\begin{tikzcd}
        & (G \ltimes \mathfrak{g}^*) \ltimes (\mathfrak{g} \ltimes \mathfrak{g}^*) \\
        G \ltimes \mathfrak{g}^* \arrow[r,"X",swap] \arrow[ur,"\Tilde{X}"] & T(G \ltimes \mathfrak{g}^*) \arrow[u,"tr_{T(G \ltimes \mathfrak{g}^*)}",swap]
\end{tikzcd}
\end{center}
\begin{align*}
    \tilde{X}(g,\mu)=(g,\mu;\mathbb{I}^{-1}(\mu),0).
\end{align*}
Finally, we put together all the ingredients assembled thus far to construct our integrator as
\begin{align}
    &t\tilde{X}(pr_{12}((\tilde{R}_d^{T^*})^{-1}(g_k,\mu_k;g_{k+1},\mu_{k+1}))=
    (\tilde{R}_d^{T^*})^{-1}(g_k,\mu_k;g_{k+1},\mu_{k+1})\nonumber\\
    &\implies t\tilde{X}(g_k,-d^{L*}_{\tau^{-1}(g_k^{-1}g_{k+1})}\tau(\mu_{k+1}))=
    (g_k,-d^{L*}_{\tau^{-1}
    (g_k^{-1}g_{k+1})}\tau(\mu_{k+1});\nonumber\\
    &\tau^{-1}(g_k^{-1}g_{k+1}),
    \mu_k-Ad^*_{g_{k+1}^{-1}g_k}(\mu_{k+1}))\nonumber\\
    &\implies (-t\mathbb{I}^{-1}(d^{L*}_{\tau^{-1}(g_k^{-1}g_{k+1})}\tau(\mu_{k+1})),0)=
    (\tau^{-1}(g_k^{-1}g_{k+1}),\mu_k\nonumber\\
    &-Ad^*_{g_{k+1}^{-1}g_k}(\mu_{k+1}))\nonumber\\
    &\implies \tau^{-1}(g_k^{-1}g_{k+1})+t\mathbb{I}^{-1}(d^{L*}_{\tau^{-1}(g_k^{-1}g_{k+1})}\tau(\mu_{k+1}))=0\nonumber\\
    &\text{and} \quad \,\mu_{k+1}=Ad^*_{g_k^{-1}g_{k+1}} \mu_k\nonumber\\
    &\implies \xi_{k,k+1}+t\mathbb{I}^{-1}(d^{L*}_{\xi_{k,k+1}}\tau(\mu_{k+1}))=0\\
    &\text{and} \quad \,\mu_{k+1}=Ad^*_{\tau(\xi_{k,k+1})} \mu_kr\\
    &\text{and} \quad \, g_{k+1}=g_k\tau(\xi_{k,k+1}).
\end{align}

\end{document}